\documentclass[oneside, a4paper,11pt,reqno]{amsart}
\usepackage{xcolor}
\usepackage{ulem}

\newtheorem{hypo1}{Hypothesis}[section]

\newtheorem{prop1}[hypo1]{Proposition}
\newtheorem{thm1}[hypo1]{Theorem}
\newtheorem{lem1}[hypo1]{Lemma}

\newtheorem{hypo}{Hypothesis}

\newtheorem{thm}[hypo]{Theorem}

\def\PP{\mathbb{P}}
\def\RR{\mathbb{R}}

\def\EE{\mathbb{E}}
\def\NN{\mathbb{N}}

\def \ind {{\bf 1}}
\newcommand{\cal}[1]{\mathcal{#1}}
\newcommand {\pare}[1] {\left( {#1} \right)}
\newcommand {\cro}[1] {\left[ {#1} \right]}
\newcommand {\va}[1] {\left| {#1} \right|}
\newcommand {\acc}[1] {\left\{ {#1} \right\}}

\newcommand {\sur}[2] { \stackrel {\scriptstyle{#1}}{#2}}
\newcommand {\ceil}[1] {\left\lceil {#1} \right\rceil}

 \title{Persistence exponent for random processes in Brownian scenery}
 \author{Fabienne Castell} 
\address{Aix-Marseille Universit\'e, CNRS, Centrale Marseille, 
Institut de Math\'ematiques de Marseille, UMR 7373, 
 13 453 Marseille. 
 France.}
\email{fabienne.castell@univ-amu.fr}

\author{Nadine Guillotin-Plantard} 
\address{Institut Camille Jordan, CNRS UMR 5208, Universit\'e de Lyon, Universit\'e Lyon 1, 43, Boulevard du 11 novembre 1918, 69622 Villeurbanne, France.}
\email{nadine.guillotin@univ-lyon1.fr}

\author{Fr\'ed\'erique Watbled} 
\address{LMBA, Universit\'e de Bretagne Sud,  Campus de Tohannic, BP 573, 56017 Vannes, France. }
\email{watbled@univ-ubs.fr}

\subjclass[2000]{60F05; 60F17; 60G15; 60G18; 60K37}
\keywords{Random scenery, First passage time, One-sided barrier problem, One-sided exit problem, Survival exponent, Persistence. \\
This research was supported by the french ANR project MEMEMO2 2010 BLAN 0125}

\medskip
 
\begin{document}

\begin{abstract}
In this paper we consider the persistence properties of random processes in Brownian scenery, which are examples of non-Markovian and
non-Gaussian processes. More precisely we study the asymptotic behaviour for 
large $T$, of the probability $\PP[ \sup_{t\in[0,T]} \Delta_t \leq 1] $
where $\Delta_t = \int_{\mathbb{R}} L_t(x) \, dW(x).$
Here $W=\acc{W(x); x\in\mathbb{R}}$ is a two-sided standard real Brownian motion and $\acc{L_t(x); x\in\mathbb{R},t\geq 0}$ is
the local time of some self-similar random process $Y$, independent from the process
$W$. We thus generalize the results of \cite{BFFN} where the increments of $Y$ were assumed to be independent.
\end{abstract}

\maketitle

\section{Introduction} 
Let $W=\acc{W(x); x\in\mathbb{R}}$ be a standard two-sided real Brownian motion and $Y=\acc{Y(t); t\ge 0}$ be a real-valued self-similar process of index $\gamma \in (0, 2)$ (i.e. for any $c >0$, $\acc{Y(c t); t\ge 0} \sur{(d)}{=} \acc{c^{\gamma} Y(t); t\ge 0}$), with stationary increments. When it exists,  
we will denote by $\acc{L_t(x); x\in\mathbb{R},t\geq 0}$ a version of the local time of the process $\acc{Y(t) ; t\geq 0}$. 
The process $L$ satisfies the occupation density formula: for each bounded measurable function $f:\RR\to\RR$ and for each $t\geq 0$,
\begin{equation}\label{odf}
\int_0^t f(Y(s))ds=\int_\RR f(x)L_t(x)dx.
\end{equation}
 The processes $W$ and $Y$ are defined on the same probability space and are assumed to be independent.
We consider {\it the random process in Brownian scenery} $\acc{\Delta_t;t\geq 0}$ defined as
\[\Delta_t = \int_{\mathbb{R}} L_t(x) \, dW(x).\]
The process $\Delta$ is itself a self-similar process of index $h$ with stationary increments, with
\[h:=  1-\frac{\gamma}{2}.\]
This process can be seen as a mixture of Gaussian processes, but it is neither Gaussian nor Markovian. It appeared independently in the mathematical literature (see \cite{KS},\cite{Borodin}), and in the physics literature (see \cite{MdM})
where it was originally introduced to model diffusion in layered white-noise velocity field. Note indeed that $\Delta_t$ can
formally be written as $\int_0^t \dot{W}(Y(s)) ds$. Thus, it  describes the horizontal motion of a tracer particle
in a (1+1)-dimensional medium, where the motion of the particle along the vertical direction $y$ is described by the process $Y$, while along the horizontal direction, the particle is driven by the white noise velocity field $\dot{W}(y)$, that depends only on the vertical coordinate $y$. The process $\acc{\Delta_t;t\geq 0}$ provides a simple example of anomalous super-diffusion, being for instance of order $t^{3/4}$ when $Y$ is a Brownian motion. In this paper, we are interested
in the {\it persistence probability} of the process $\Delta$, i.e. the asymptotic
behaviour of 
$${\mathcal F}(T): = \PP \Big[ \sup_{t\in[0,T]} \Delta_t \leq 1\Big] $$
as $T\rightarrow +\infty$. The study of the one-sided supremum process of random processes is a classical issue
in probability. One usually gets polynomial decay of the persistence probability: $\cal F(T) \asymp T^{-\theta}$ for 
a non-negative $\theta$, often called {\it persistence exponent}, or {\it survival exponent}. Classical examples,
where this exponent can be computed, include random walks or L\'evy processes, and we refer the reader to the recent survey papers \cite{AS, BMS} for an account of models where the persistence exponent is known.  However, there are relevant 
physical situations where this exponent remains unknown (see for instance \cite{Maj1,BMS}). The persistence exponent of the Brownian motion in Brownian scenery was studied by Redner \cite{Red1,Red2}, 
and Majumdar \cite{Maj}. Based on physical
arguments, numerical simulations and analogy with fractional Brownian motion, Redner and Majumdar conjectured the
value of the persistence exponent. In \cite{BFFN}, their conjecture was proved up to logarithmic factors, when the
process $Y$ is a stable L\'evy process with index $\delta \in (1,2]$. The proof of \cite{BFFN} depends heavily on the increments independence
of the process $Y$ and the question raises if it is possible to compute the persistence exponent without it. 
The aim of this paper is to answer this question, and to provide  assumptions on $Y$ allowing to compute the persistence exponent of the
random process in Brownian scenery. 
\begin{description}
\item[(H1)] There exists a continuous version $\acc{L_t(x); x\in\mathbb{R},t\geq 0}$ of the local time of $Y$;
\item[(H2)] $Y$ is a self-similar process of index $\gamma \in (0, 2)$;
\item[(H3)] For every $T>0$, $\{Y(T-t) - Y(T); t\in [0,T]\}$ has the same distribution as either $\{Y(t); t\in[0,T]\}$ or $\{-Y(t); t\in [0,T]\}$;
\item[(H4)] Let $V_1:=\int L_1^2(x) \, dx$ be the self-intersection local time of $Y$. There exist a real number 
$\alpha > 1$, and positive constants $C,c$ such that for any $x \ge 0$,
$$\PP\cro{V_1 \ge x} \le C \exp(-cx^{\alpha}) \, ;
$$
\item[(H5)]  There exist a real number 
$\beta > 0$, and positive constants $C,c$ such that for any $x > 0$,
$$\PP\cro{V_1 \le x} \le C \exp(-cx^{-\beta}) \, .
$$
\end{description}
Our main result is the following one.
\begin{thm}
\label{theoMk}
Assume (H1) to  (H5)  hold. There exists a constant $c>0$,  such that for large enough $T$ , 
\begin{equation}
\label{eqMk}
 T^{-\gamma/2} (\ln T)^{-c}\leq \PP\Big[ \sup_{t\in[0,T]} \Delta_t\leq 1\Big]
\leq  T^{-\gamma/2} (\ln T)^{+c}.
\end{equation}
\end{thm}

 Our main achievement in this paper is to remove the  increments independence assumption on 
$\{Y(t); t \geq 0\}$, which was crucial in \cite{BFFN}. 
Another interesting issue would be to obtain the persistence
exponent for $\{W(x); x \in \RR\}$ being a stable L\'evy process. This seems delicate,  
our proof  relying on the fact that the process  $\{\Delta(t); t \geq 0\}$ is a Gaussian process, conditionally 
to $\{Y(t); t \geq 0\}$.

The paper is organized as follows. Section \ref{exemples} provides three examples of processes satisfying 
(H1) to (H5), including the stable L\'evy process with index $\delta \in (1,2]$,
and thus generalizing the result of \cite{BFFN}.
Section \ref{aux} states some useful properties of the process $\Delta$. 
Section \ref{Pf} is devoted to the proof of Theorem \ref{theoMk}.

\section{Some examples}
\label{exemples}
\subsection{A sufficient condition for (H5)}
In this section we assume that (H1) holds and we provide a sufficient condition on $\{Y(t); t \geq 0\}$ allowing to check (H5). 
Let us first compare $V_1$ with $\displaystyle\max_{s\in [0,1]}|Y(s)|$.
By noting that
$$1=\int_{\{|x|\leq \max_{s\in [0,1]}|Y(s)|\}} L_1(x) dx\leq V_{1}^{1/2} \sqrt{2\max_{s\in [0,1]}|Y(s)|}  ,$$
we deduce that
\begin{equation}\label{comparaison}
\frac{1}{2\displaystyle\max_{s \in [0,1]} \va{Y(s)}}\leq V_1.
\end{equation}
 
\begin{lem1}\label{CS}
Assume that (H1) holds and that there exist positive constants $C$ and $c$, and a real number $\beta >0$ such that for all $x>0$, 
$$ \PP\cro{\max_{s \in [0,1]} \va{Y(s)} \ge x} \leq C \exp\pare{-c x^{\beta}} \,.$$
Then (H5) holds true.
\end{lem1}
\begin{proof}
It is a direct consequence of \eqref{comparaison} that for all $x>0$,
$$\PP\cro{V_1 \le x} \leq \PP\cro{\max_{s \in [0,1]} \va{Y(s)} \ge \frac{1}{2x}} 
\leq C \exp(-c2^{-\beta} x^{-\beta}) \, . 
$$
 \end{proof} 
 
\subsection{Stable L\'evy processes.}
A process $Y=\acc{Y(t); t\ge 0}$ is  a strictly stable L\'evy process with index $\delta\in (1,2]$ if it is a 
process with stationary and independent increments,  such that $Y(0)=0$ and for any $t \geq 0$ and 
$u \in \RR$, 
\begin{equation}
\label{cdfY}
 \mathbb{E}[e^{iuY(t)}]=\exp\Big\{- c |u|^{\delta} t  \Big(1 + i \zeta \   \mbox{\rm sgn} (u) \tan(\pi \delta /2) \Big)\Big\} \, ,
\end{equation}
where $ \zeta \in [-1,1]$, and $c$ is a positive scale parameter. 
\begin{lem1}
The stable L\'evy process with index  $\delta\in (1,2]$ satisfies (H1) to (H5) with $\gamma=\frac{1}{\delta}$,
$\alpha=\delta$, and $\beta=\frac{\delta}{2\delta-1}$.
\end{lem1}
\begin{proof}
It is immediate that the L\'evy process is self-similar of index $1/\delta$ and that for every $T>0$, 
$$\{Y(T-t) - Y(T); t\in [0,T]\} \overset{\text(d)}{=}\{-Y(t); t\in [0,T]\}$$
then (H2) and (H3) are satisfied.\\*
A  continuous version of the local time exists, since $\delta > 1$ (this was proved by Boylan \cite{Boylan}), so it satisfies (H1).
Moreover, Corollary 5.6 in \cite{KL} states that 
there exist positive constants $C$ and $\xi$ s.t. for every $x>0$,
$$\PP[V_1\geq x]  \leq  C e^{-\xi x^{\delta} },$$
which gives (H4).
To prove (H5) we let $V_t:=\int L_t^2(x) \, dx$ be the self-intersection local time of $Y$ up to time $t$ for $t\geq 0$, and we show that there exists $l>0$ such that
\begin{equation}\label{H5levy}
\PP\cro{V_t\leq 1}\leq e^{-lt} \textrm{ for } t\textrm{ large enough}.
\end{equation}
As $\{V_t; t\geq 0\}$ is self-similar of index $\frac{2\delta-1}{\delta}$,
(\ref{H5levy}) implies that for every positive $\varepsilon$ small enough, 
$$\PP[V_1\leq \varepsilon]=\PP[V_{\varepsilon^{-\frac{\delta}{2\delta-1}}}\leq 1]\leq e^{-l\varepsilon^{-\frac{\delta}{2\delta-1}}},$$
which gives (H5).
To prove \eqref{H5levy} we show that the function defined on $\RR^+$ by 
$$f(t)=\log \PP[V_t\leq 1]$$
is subadditive. Let us fix $s$, $t$ in $[0,+\infty)$.
We consider the process $Y^{(s)}:=\{Y^{(s)}_u; u\geq 0\}$ defined by
$$Y^{(s)}_u=Y(u+s)-Y(s),$$
its local time $\{L^{(s)}_u(x);x\in\RR,u\geq 0\}$, and its self-intersection
local time $\{V^{(s)}_u; u\geq 0\}$.
Note that
$$L_{t+s}(x)=L_s(x)+L^{(s)}_t(x-Y(s)).$$
Hence
$$V_{t+s}\geq V_s+V^{(s)}_t.$$
The process $Y$ being a L\'evy process, $V^{(s)}_t$ and $V_s$ are independent
and $V^{(s)}_t$ has the same law as $V_t$. Therefore
$$\PP[V_{t+s}\leq 1]\leq \PP[V_s\leq 1; V^{(s)}_t\leq 1]=
\PP[V_s\leq 1]\PP[V_t\leq 1].$$
Thus $f$ is subadditive, which implies that $\frac{f(t)}{t}$ converges, as $t\rightarrow +\infty$, towards a limit $-l=\inf_{t>0}\frac{f(t)}{t}$ (see \cite{Hamm}). 
It remains to show that $l$ is strictly positive. As $-l$ is less than or equal to $f(1)$,
it is enough to show that $\PP[V_1>1]>0$.
But thanks to \eqref{comparaison} and Proposition 10.3 of \cite{Fristedt},
$$\PP[V_1>1]\geq \PP\cro{\max_{s \in [0,1]} \va{Y(s)}<\frac{1}{2}}>0,$$
and this concludes the proof of \eqref{H5levy}.
\end{proof}

\subsection{Fractional Brownian motion.}
The fractional Brownian motion of Hurst index $H \in (0,1)$ is the real centered Gaussian process 
$\acc{B_H(t); t \geq 0}$ 
with covariance function 
$$\EE[B_H(t)B_H(s)]=  \frac{1}{2} (t^{2H}+s^{2H} -\va{t-s}^{2H}) \, .
$$
\begin{lem1}
The fractional Brownian motion with Hurst index $H \in (0,1)$ satisfies (H1) to (H5) with $\gamma=H$,
 $\alpha=1/H$, $\beta =2$.
\end{lem1}
\begin{proof}
It follows readily from the definition that $B_H$ is self-similar with index $H$ and satisfies $(H3)$. 
The existence of a jointly continuous version of its local time process for $H\in (0,1)$ is a classical fact
(see for instance paragraph 22 in \cite{GH}).  
Theorem 1 of \cite{HNS} asserts that there exists 
$\lambda_0 > 0 $ such that for every $0< \lambda < \lambda_0$,
\begin{equation}
\label{momentexpV_1}
\EE\cro{e^{\lambda V_1^{1/H}}} < \infty \, ,
\end{equation}
which implies (H4) with $\alpha =1/H$. 
Finally,  (H5) follows from Lemma \ref{CS} and Fernique's estimation (\cite{Fernique}, Theorem 4.1.1): there exists
$c_F>0$ such that for any $x\geq \sqrt 5$, 
\begin{equation}
\label{ineg.Fernique}
\PP[\max_{s \in [0,1]} \va{B_H(s)} \geq c_Fx]\leq 10\int_x^{+\infty}e^{-\frac{v^2}{2}}dv \, .
\end{equation}
This implies (H5) with $\beta =2$. 
\end{proof}

\subsection{Iterated Brownian motion.}
Let $\acc{B(x); x\in \RR}$ be a two-sided real standard Brownian motion, and let $\acc{\tilde{B}(t); t \geq 0}$ be
a real
standard Brownian motion, independent of $\acc{B(x); x\in \RR}$. The iterated Brownian motion is the process
$\acc{Y(t); t \geq 0}$ defined by 
$$ Y(t) = B(\tilde{B}(t)) \, .
$$
\begin{lem1}
The iterated Brownian motion satisfies (H1) to (H5) with $\gamma=1/4$,
and $\alpha=\beta=4/3$.
\end{lem1}
\begin{proof}
The self-similarity of the iterated Brownian motion is a direct consequence 
of the self-similarity and independence of both Brownian motions. The assertion $(H3)$ follows once again from the independence of $B$ and $\tilde{B}$, the increments stationarity of 
$B$ and the fact that $\tilde{B}$ satisfies $(H3)$.

 The existence and joint continuity
of the local times of iterated Brownian motion were proved in \cite{BK,CCFR}.  Let us prove 
(H5). For $x > 0$, 
\begin{eqnarray}
\nonumber \PP\cro{\max_{s\in [0,1]} \va{Y(s)} \ge x} 
 & \leq &  \PP\cro{\max_{s\in [0,1]} \va{\tilde{B}(s)} \ge x^{2/3}} 
  + \PP\cro{\max_{\va{u} \leq x^{2/3}}  \va{B(u)} \ge x}  
 \\
 \nonumber &  \leq  & \PP\cro{\max_{s\in [0,1]} \va{\tilde{B}(s)}  \ge x^{2/3}} 
  + \PP\cro{\max_{u \in [0, x^{2/3}]}  \va{B(u)} \ge x} 
\\
\nonumber &   \leq  & 2 \,  \PP\cro{\max_{s\in [0,1]} \va{B(s)}  \ge x^{2/3}}  
 \\
 &  \leq  & 4 \,  \PP\cro{\max_{s\in [0,1]} B(s)  \ge x^{2/3}} 
  = 4 \, \PP\cro{\va{B(1)}  \ge x^{2/3}}
   \, .  
   \label{IBMMax}
\end{eqnarray}
This proves (H5) with $\beta=4/3$ using Lemma \ref{CS}.

 To prove (H4), we use the uniform norm on the local times of iterated Brownian motion proved in Lemma 4
  of \cite{Xiao}: there exists a constant $K > 0$ such that for any even integer $n$,  
 \begin{equation}
 \label{xiao.UB}
 \sup_{x \in \RR}  \EE \cro{(L_1(x))^n} \leq K^n (n!)^{3/4} \, . 
\end{equation}
Since $\int_{\RR} L_1(x) \, dx = 1$, H\"older's inequality implies  that for any integer $n$, 
$$V_1^n = \pare{\int_{\RR} L_1(x)^2 \, dx }^n \leq \int_{\RR} L_1(x)^{n+1}  \, dx = 
\int_{-\max_{s \in [0,1]} \va{Y(s)}} ^{\max_{s \in [0,1]}  \va{Y(s)}} \, L_1(x)^{n+1}  \, dx \, . 
$$
Therefore, 
\begin{eqnarray}
\nonumber
\EE\cro{V_1^n} & \leq  & \int_{\RR} \EE\cro{ L_1(x)^{n+1} \ind_{\va{x} \leq \max_{s \in [0,1]}  \va{Y(s)} }} \, dx
\\
\nonumber
& \leq & \int_{\RR}  \sqrt{ \EE\cro{ L_1(x)^{2(n+1)}}} \sqrt{\PP\cro{\max_{s \in [0,1]}  \va{Y(s)} \geq \va{x}}} \, dx 
\\
& \leq & K^{n+1}  \, ( (2n+2)!)^{3/8} \int_{\RR} \sqrt{\PP\cro{\max_{s \in [0,1]}  \va{Y(s)} \geq \va{x}}} \, dx \, , 
\label{IBMSILT}
\end{eqnarray}
where we used \eqref{xiao.UB} in the last inequality.  By \eqref{IBMMax}, the integral in  \eqref{IBMSILT} is finite. We deduce then from Stirling's formula that there
exists a constant $C> 0$ such that for any integer $n$, 
$$ \EE\cro{V_1^n}  \leq C^n  n^{\frac{3}{4} n} \, .
$$
Hence, for any $x> 0$, and any integer $n$, 
$$\PP\cro{V_1 \geq x} \leq x^{-n} C^n  n^{\frac{3}{4} n} \, .
$$ 
Optimizing over the values of $n$  leads to take $n = \ceil{e^{-1} \pare{\frac{x}{C}}^{4/3}}$, so that 
$ x^{-n} C^n  n^{\frac{3}{4} n}  \simeq \exp(-\frac{3}{4e} \pare{\frac{x}{C}}^{4/3})$ for large $x$. 
This proves that the iterated Brownian motion satisfies (H4) with $\alpha=4/3$.

\end{proof}

\section{Auxiliary statements on $(\Delta_t; t \geq 0)$}
\label{aux}

For a certain class of stochastic processes $\acc{X_t; t\geq 0}$ (to be specified below),
 Molchan \cite{Molchan1999} proved that the asymptotic behavior of 
$$\PP\Big[ \sup_{t\in[0,T]} X_t \leq 1\Big] $$  is related to the quantity 
$$I(T) :=  \EE\left[\left( \int_0^T e^{X_t}\ dt\right)^{-1} \right].$$
We refer to \cite{Aurzada} where the relationship between both quantities is clearly explained as well as the heuristics. Moreover $I(T)$ can be estimated 
by the following result.
\begin{thm1}[Statement 1, \cite{Molchan1999}]\label{Molchan}
Let $\acc{X_t; t\geq 0}$ be a continuous process, self-similar with index $h > 0$ such that for every $T>0$, 
$$\{X_{T-t} - X_{T}; t\in [0,T]\} \overset{\text(d)}{=}\{X_t; t\in [0,T]\}.$$
Moreover assume that for every $\theta>0$,
$$\EE\Big[\exp{\Big(\theta \max_{t\in[0,1]} |X_t|\Big)} \Big] <+\infty.$$
Then, as $T\rightarrow +\infty$,
$$ \EE\left[\left( \int_0^T e^{X_t}\ dt\right)^{-1} \right] = h T^{-(1-h)} \Big( \EE\Big[\max_{t\in[0,1]} X_t \Big] +o(1) \Big).$$
\end{thm1}
\noindent Before applying this result to our random process $\Delta$ we first establish some useful facts concerning it. First we show that the process satisfies
$$\{\Delta_{T-t} - \Delta_{T}; t\in [0,T]\} \overset{\text(d)}{=}\{\Delta_t; t\in [0,T]\}$$ 
 and has stationary increments (Lemma \ref{statincr}).
Next we provide an exponential upper bound for the tail of $\Delta_1$,  from which we deduce the Kolmogorov-Centsov continuity criterion (Lemma \ref{PropDelta1}).
Finally we show that the process satisfies a maximal inequality (Lemma \ref{PropInegmax}).

\begin{lem1}\label{statincr}
Assume $(H1)$ and $(H3)$. The process $\Delta$ satisfies for every $T>0$, 
 \begin{equation}
\label{RenvTemps}
 \{\Delta_{T-t} - \Delta_{T}; t\in [0,T]\} \overset{\text(d)}{=}\{\Delta_t; t\in [0,T]\} \, . 
\end{equation}
Therefore, it  has stationary increments, i.e. for every $s>0$, $\{\Delta_{t+s} - \Delta_{s}; t\geq 0\} \overset{\text(d)}{=}\{\Delta_t; t\geq 0\}$.
\end{lem1}
\begin{proof}
Conditionally to $\{Y(t); t \geq 0\}$, the process $\{\Delta_{T-t} - \Delta_{T}; t\in [0,T]\}$ is a centered 
Gaussian process with (random) covariance function:
$$ C_T(s,t) = \int_{\RR} (L_T(x)-L_{T-t}(x)) \, (L_T(x)-L_{T-s}(x)) \, dx  \, .
$$ 
Note that $L_T(x)-L_{T-t}(x)= L^{(T)}_t(x-Y(T))$, where $\{L^{(T)}_t(x); t \in [0,T], x \in \RR\}$ is the local 
time process of $Y^{(T)}_t :=Y(T-t)-Y(T) $. Hence, 
$ C_T(s,t)= \int_{\RR} L^{(T)}_t(x) L^{(T)}_s(x) \, dx$. Now,  if $Y^{(T)}$ has the same distribution as 
$Y$  (respectively $-Y$), then $\{L^{(T)}_t(x); t \in [0,T],  x \in \RR\}$  has the same law as $\{L_t(x); t \in [0,T],  x \in \RR \}$ (respectively $\{L_t(-x); t \in [0,T],  x \in \RR\}$).  We deduce then that in both cases, 
$\{C_T(s,t); s, t \in [0,T]\}$ is distributed as the conditional covariance of  $\{\Delta_{t} ; t\in [0,T]\}$ with 
respect to $Y$. Hence, 
 $$\{\Delta_{T-t} - \Delta_{T}; t\in [0,T]\} \overset{\text(d)}{=}\{\Delta_t; t\in [0,T]\} \, .$$ 
Concerning the increments stationarity, fix $s >0$, and $0 < t_1 < \cdots < t_n$. Let $T > t_n + s$ be fixed. 
$$
\begin{aligned}
& (\Delta_{t_1+s} - \Delta_s, \cdots, \Delta_{t_n+s} - \Delta_s)
\\
& \hspace*{1cm} \overset{\text(d)}{=} 
(\Delta_{T-t_1-s} - \Delta_{T-s}, \cdots, \Delta_{T-t_n-s} - \Delta_{T-s})
\mbox{ by time reversal at time } T \, ,
\\
&   \hspace*{1cm} \overset{\text(d)}{=} 
(\Delta_{t_1} , \cdots, \Delta_{t_n})
\mbox{ by time reversal at time } T-s \, .
\end{aligned}
$$
 
\end{proof}

\begin{lem1}\label{PropDelta1}
Assume (H1) and (H4) hold. There exist $C>0$ and $\delta>0$ such that for any $x>0$,
\begin{equation}\label{queueDelta_1}\PP[\Delta_1\geq x]\leq C\exp(-\delta x^{\frac{2\alpha}{1+\alpha}}).
\end{equation}
Then, for every $a\geq 1$,
\begin{equation}\label{KC}
(\EE | \Delta_t -\Delta_s |^a )^{1/a}   = C(a) |t-s|^{h}, \    t,s \geq 0
\end{equation}
where $C(a)\leq c  a^{\nu}$ with $\nu := \frac{1}{2} \pare{1+\frac{1}{\alpha}}$. 
In particular the process satisfies the  Kolmogorov-Centsov criterion of continuity.
\end{lem1}
\begin{proof}
Conditionally to the process $Y$, the random variable $\Delta_1$ is a real centered
Gaussian variable with variance $V_1.$
For each $u\in \RR$, let
$$\Phi(u)=\frac{1}{\sqrt{2\pi}}\int_u^{+\infty}e^{-\frac{s^2}{2}}ds.$$
Then for $x>0$ and $\theta>0$,
$$\begin{aligned}
\PP[\Delta_1\geq x]=\int_0^{+\infty}\Phi(xz^{-1/2})\PP_{V_1}(dz)
&=\int_0^{x^\theta}\Phi(xz^{-1/2})\PP_{V_1}(dz)
+\int_{x^\theta}^{+\infty}\Phi(xz^{-1/2})\PP_{V_1}(dz)\\
&\leq \Phi(x^{1-\frac{\theta}{2}})+\PP[V_1\geq x^\theta],
\end{aligned}$$
where we used that the function $z\mapsto \Phi(xz^{-1/2})$ is non decreasing with values in $[0,1]$.
Using (H4) and  the classical inequality
$$\Phi(u)\leq \frac{1}{\sqrt{2\pi}u}e^ {-u^2/2}\textrm{ for every }u>0,$$
and choosing $\theta=\frac{2}{1+\alpha}$,
we obtain that for $x$  large enough,
$$\PP[\Delta_1\geq x]\leq \frac{x^{-\alpha/(1+\alpha)}}{\sqrt{2\pi}}e^{-\frac{1}{2}x^{2\alpha/(1+\alpha)}}
+C e^ {-cx^{2\alpha/(1+\alpha)}}\leq  C e^ {-\delta x^{2\alpha/(1+\alpha)}},$$
with $\delta=\min(c,1/2)$. This implies \eqref{queueDelta_1} for every $x>0$.

Let us prove (\ref{KC}). The increments of the process $\Delta$ being stationary (see Lemma (\ref{statincr})), by self-similarity, we have for every $t,s\geq 0$,
\begin{eqnarray*}
\EE[ |\Delta_t -\Delta_s|^a]&=& |t-s|^{h a}\   \EE[ |\Delta_1|^a] \, . \\ 
\end{eqnarray*}
 Now, using \eqref{queueDelta_1} and the symmetry of $\Delta_1$,
\begin{eqnarray}
\nonumber
\EE[ |\Delta_1|^a]
&= & \int_0^{+\infty} \PP[ |\Delta_1| \geq x^{1/a}] \, dx 
 \leq  2  \int_0^{+\infty} \PP[ \Delta_1 \geq x^{1/a}] \, dx
\\
\label{momentDelta}
& \leq &C \int_0^{+\infty} \exp(-\delta x^{\frac{2\alpha}{a(1+\alpha)}})  \, dx
= C a \delta^{-\frac{a(1+\alpha)}{2\alpha} } \Gamma\pare{\frac{a(1+\alpha)}{2\alpha}} \, . 
\end{eqnarray}
Hence, it follows from Stirling's formula that $C(a)= \EE[ |\Delta_1|^a]^{1/a}  \leq C a^{(1+\alpha)/(2\alpha)}$. 

\end{proof}

\begin{lem1}\label{PropInegmax}
Assume (H1) and (H4) hold. Let $T,x\geq 0$.
Then
\begin{equation}\label{inegmaxcond}
\PP[\max_{s\in [0,T]}\Delta_s\geq x|Y]\leq 2 \PP[\Delta_T\geq x|Y],
\end{equation}
\begin{equation}\label{inegmax}
\PP[\max_{s\in [0,T]}\Delta_s\geq x]\leq 2 \PP[\Delta_T\geq x].
\end{equation}
\end{lem1}
\begin{proof}
Conditionally to the process $Y$, the process $\{\Delta_t; t\geq 0\}$ is a centered
Gaussian process on $\RR$ with covariance function
$$\EE\cro{\Delta_s \Delta_t| Y} =\int_\RR L_s(x)L_t(x)dx.$$
Moreover for any $t\geq s\geq 0$,
$$\begin{aligned}
\EE\Big[\Delta_t^2|Y\Big]-\EE\Big[\Delta_s^2|Y\Big]-
\EE\Big[(\Delta_t-\Delta_s)^2|Y\Big]
&=\int_\RR (L_t^2(x)-L_s^2(x)-(L_t(x)-L_s(x))^2)dx\\
&=\int_\RR 2L_s(x)(L_t(x)-L_s(x))dx\geq 0,
\end{aligned}$$
hence applying Proposition 2.2 in \cite{KL}, we deduce the inequality \eqref{inegmaxcond}.
By integrating we obtain the maximal inequality \eqref{inegmax}.
\end{proof}

We are now in position to use Theorem \ref{Molchan} and to state the main result of this section.
\begin{prop1}\label{Prop1}
Assume (H1) to (H4) hold. As $T\rightarrow +\infty$,
$$ \EE\left[\left( \int_0^T e^{\Delta_t}\ dt\right)^{-1} \right] =\left(1-\frac{\gamma}{2}\right) T^{-\gamma/2} 
\Big( \EE\Big[\max_{t\in[0,1]} \Delta_t \Big] +o(1) \Big).$$
\end{prop1}
\begin{proof}
It follows easily from assumption (H2), that the process $\acc{\Delta_t; t\geq 0}$ is self-similar with index $h:=1-\frac{\gamma}{2}$. From Lemma \ref{statincr}, it satisfies \eqref{RenvTemps}. From \eqref{KC} and Kolmogorov-Centsov continuity criterion, we can assume that it is continuous. Hence, to apply Theorem \ref{Molchan}, it is enough to prove that for every $\theta>0$,

\begin{equation}\label{finito}
\EE\Big[\exp{\Big(\theta \max_{t\in[0,1]} |\Delta_t|\Big)} \Big] <+\infty.
\end{equation}
Let $\theta>0$. We have\\*
\begin{eqnarray*}
\EE\Big[\exp{\Big(\theta \max_{t\in[0,1]} |\Delta_t|\Big)} \Big]  &=& \int_0^{\infty} \PP\Big[ \exp{(\theta \max_{t\in[0,1]} |\Delta_t|)} \geq \lambda  \Big] d\lambda\\
&\leq & 2+ \int_2^{\infty} \PP\Big[ \max_{t\in[0,1]} |\Delta_t| \geq \frac{\ln (\lambda)}{\theta}  \Big] d\lambda.
\end{eqnarray*}
Since the process $\acc{\Delta_t; t \geq 0}$ is symmetric and satisfies the maximal inequality \eqref{inegmax} of Lemma \ref{PropInegmax},
\begin{eqnarray*} 
\PP\Big[ \max_{t\in[0,1]} |\Delta_t| \geq \frac{\ln (\lambda)}{\theta}  \Big] &\leq & 2\  \PP\Big[ \max_{t\in[0,1]} \Delta_t \geq \frac{\ln (\lambda)}{\theta}  \Big]\\
&\leq & 4\  \PP[ \Delta_1 \geq (\ln \lambda)/\theta].
\end{eqnarray*} 
We apply the inequality \eqref{queueDelta_1} of Lemma \ref{PropDelta1}, and
since the function $\lambda\rightarrow  \exp( - \delta ((\ln \lambda)/\theta)^{2\alpha/ (1+\alpha)} ) $ 
is integrable at infinity for any $\alpha > 1$, the assertion (\ref{finito}) follows.
\end{proof}

\section{Proof of Theorem \ref{theoMk}} 
\label{Pf}
 In this section, we prove upper and lower bounds on the persistence probability 
$\PP\cro{\sup_{t\in[0,T)} \Delta(t) \leq 1}$. 

In the case of fractional Brownian motion $\{ B_{H} (t); t \geq 0\}$,
Aurzada's proof of the lower bound  (see \cite{Aurzada}),  rests on both following arguments: the fractional Brownian motion satisfies the hypothesis of Theorem \ref{Molchan}, and it satisfies the equality  (valid for $a$ large enough)
\begin{equation} \label{ing}
(\EE | B_{H} (t) -B_{H}(s) |^a )^{1/a}   = C(a) |t-s|^H, \    t,s \geq 0
\end{equation}
with $C(a) \leq c a^{\nu}$, for some $c$ and $\nu>0$. 
We showed in Lemmas \ref{statincr}  and \ref{PropDelta1} that these conditions are also satisfied 
by our random process $\Delta$. 
Therefore the proof of Aurzada \cite{Aurzada} allows us to derive the lower bound in Theorem \ref{theoMk}. 


As in \cite{Molchan1999} and \cite{Aurzada}, the main idea in the proof of the upper
bound in \eqref{eqMk}, is to bound $I(T)$  from below 
by restricting the expectation  to a well-chosen set of paths. 

\noindent Conditionally to $Y$, the process $\acc{\Delta_t; t \ge 0}$
is a centered Gaussian process such that for every $0\le t<s$, 
$$\EE[\Delta_t \Delta_s|  Y] = \int_{\mathbb R} L_t(x) L_s(x)\, dx  \geq 0 , $$
$$ \EE[\Delta_t (\Delta_s - \Delta_t)|  Y] = \int_{\mathbb R} L_t(x) ( L_s(x) - L_t(x)) \, dx  \geq 0 ,$$
since $t \rightarrow L_t(x)$ is a.s. increasing for all $x \in {\mathbb R}$. 
It follows then from 
Slepian's lemma,  that for every $0\le u<v<w$ and every real numbers $a,b$, 
\begin{equation}\label{slepian1}
\PP\left[\sup_{t\in[u,v]}\Delta_t\le a,\ \sup_{t\in[v,w]}\Delta_t\le b\Big|  Y \right]\ge  
    \PP\left[\sup_{t\in[u,v]}\Delta_t\le a\Big|  Y \right]
    \PP\left[\sup_{t\in[v,w]}\Delta_t\le b\Big|  Y \right]
\end{equation}
\begin{equation}\label{slepian2}
\PP\left[\sup_{t\in[u,v]}\Delta_t\le a,\ \sup_{t\in[v,w]}(\Delta_t-\Delta_v)\le b \Big|  Y \right]\ge  
    \PP\left[\sup_{t\in[u,v]}\Delta_t\le a \Big|  Y \right]\PP\left[\sup_{t\in[v,w]}(\Delta_t-\Delta_v)\le b\Big|  Y\right].
\end{equation}

Let $a_T=(\ln T)^a$, where $a>0$ will be chosen later; let 
$\alpha_T:=\frac{\ln T}{1+\ln T}$, which belongs to $(0,1)$ and 
set $\beta_T=\sqrt{V_{a_T}}$ where $V_{a_T}:=\int_\RR L_{a_T}^2(x)dx$.
Let us define the random function 
$$\phi(t):= \left\{\begin{array}{ll}
  1  & \mbox{ for } 0 \leq t < a_T \, , \\
1-\beta_T  &  \mbox{ for } a_T\leq t\leq T \, ,
\end{array} 
\right.
$$ 
which is $Y$-measurable.
Clearly, we have 
$$ \EE\left[\left( \int_0^T e^{\Delta_t}\ dt\right)^{-\alpha_T} \Big| Y \right] \geq 
\left( \int_0^T e^{\phi(t)}\ dt\right)^{-\alpha_T} \PP\Big[\forall t\in [0,T], \Delta_t\leq \phi(t)\Big|  Y \Big].$$
By Slepian's lemma (see (\ref{slepian1})), we have
\[
 \PP\Big[\forall t\in [0,T], \Delta_t\leq \phi(t)\Big|  Y \Big]
 \geq   \PP\Big[\forall t\in [0,a_T], \Delta_t\leq 1\Big|  Y  \Big] \,  \PP\Big[\forall t\in [a_T,T], 
\Delta_t\leq 1-\beta_T\Big|  Y \Big].
\]
Remark that
\begin{eqnarray*}
\PP\Big[\forall t\in [a_T,T], \Delta_t\leq 1-\beta_T \Big|  Y \Big] 
& \geq & \PP\Big[\Delta_{a_T} 
\leq -\beta_T ;\forall t\in [a_T,T], \Delta_t- \Delta_{a_T} \leq 1\Big|  Y \Big]
\\
& \geq & \PP\Big[\Delta_{a_T} \leq -\beta_T\Big| Y \Big]  \PP\Big[\forall t\in [a_T,T], 
\Delta_t- \Delta_{a_T} \leq 1\Big|  Y \Big],
\end{eqnarray*}
by Slepian's lemma (see (\ref{slepian2})). Note that 
$$ \PP[\Delta_{a_T} \leq -\beta_T| Y] = \Phi(\beta_T V_{a_T}^{-1/2})= \Phi(1).$$

Moreover, it is easy to check that when $T$ goes to infinity,
$$\int_0^T e^{\phi(t)}\ dt = O (a_T+Te^{-\beta_T}).$$
In the following $C$ is a constant whose value may change but does not depend on $T$. 
Then we can write that for $T$ large enough
\begin{equation}\label{firststep}
\EE\left[\left(\int_0^T e^{\Delta_t}dt\right)^{-\alpha_T}\Big| Y \right]
\geq C(a_T+Te^ {-\beta_T})^{-\alpha_T} \, 
\PP\big[\sup_{t\in[0,a_T]}\Delta_t\leq 1\big| Y \big]\,
\PP\big[\sup_{t\in[a_T,T]}(\Delta_t-\Delta_{a_T})\leq 1 \big| Y \big].
\end{equation}
Next we use the maximal inequality \eqref{inegmaxcond} of Lemma \ref{PropInegmax}
to write
$$\PP\big[\sup_{t\in[0,a_T]}\Delta_t\leq 1\big| Y\big]
=1-\PP\big[\sup_{t\in[0,a_T]}\Delta_t\geq 1\big| Y \big]
\geq 1-2\PP\big[\Delta_{a_T}\geq 1\big| Y \big]
=\PP[|Z|\leq V_{a_T}^{-1/2}| Y]$$
where $Z$ is a Gaussian variable $\mathcal{N}(0,1)$ independent of $Y$, from which we deduce that there
exists a constant $c>0$ such that 
\begin{equation}\label{morceau1}
\PP\big[\sup_{t\in[0,a_T]}\Delta_t\leq 1\big| Y \big]
\geq c \min( V_{a_T}^{-1/2} ,1).
\end{equation} 
Injecting \eqref{morceau1} into \eqref{firststep}
we get that for $T$ large enough,
$$\PP\big[\sup_{t\in[a_T,T]}(\Delta_t-\Delta_{a_T})\leq 1 \big|Y \big]
\leq
C\EE\left[\left. \left(\int_0^T e^{\Delta_t}dt\right)^{-\alpha_T}\right| Y \right]
(a_T+Te^ {-V_{a_T}^{1/2}})^{\alpha_T}\,\max( V_{a_T}^{1/2},1) .$$
By integrating and using successively H\"{o}lder's inequality with $p_T=\frac{1}{\alpha_T}$,
$\frac{1}{q_T}+\frac{1}{p_T}=1$, Jensen's inequality, the inequality
$(x+y)^{\alpha_T}\leq x^{\alpha_T}+y^{\alpha_T}$ for $x,y>0$, and Proposition \ref{Prop1},
we get that for $T$ large enough,
$$\begin{aligned}
\PP\big[\sup_{t\in[a_T,T]}&(\Delta_t-\Delta_{a_T})\leq 1\big]\\
&\leq
C \EE \cro{\pare{\int_0^T e^{\Delta_t}dt}^{-1}}^{1/p_T}
\EE\left[\pare{a_T^{\alpha_T} \max(V_{a_T}^{1/2},1) +T^{\alpha_T}
e^ {-\alpha_T V_{a_T}^{1/2}} \max(V_{a_T}^{1/2},1) }^{q_T}\right]^{1/q_T}
\\
&\leq CT^{-\frac{\gamma}{2p_T}}\|f_1+f_2\|_{q_T}\leq CT^{-\frac{\gamma}{2p_T}}(\|f_1\|_{q_T}+\|f_2\|_{q_T}),
\end{aligned}$$
with $$f_1=a_T^{\alpha_T} \max(V_{a_T}^{1/2},1) ,\quad f_2=T^{\alpha_T}
e^ {-\alpha_T V_{a_T}^{1/2}}\max(V_{a_T}^{1/2},1) .$$
The lefthand term is greater than the quantity we want to bound from above, since by stationarity,
$$\PP\big[\sup_{t\in[a_T,T]}(\Delta_t-\Delta_{a_T})\leq 1\big]
=\PP\big[\sup_{t\in[0,T-a_T]}\Delta_ {t}\leq 1\big]
\geq \PP\big[\sup_{t\in[0,T]}\Delta_ {t}\leq 1\big].$$
Concerning the righthand term, we recall that $\alpha_T=\frac{1}{p_T}=\frac{\ln T}{1+\ln T}$ and 
$\frac{1}{q_T}=1-\alpha_T=\frac{1}{1+\ln T}$. Hence, when $T$ goes to infinity, 
$T^{-\frac{\gamma}{2p_T}} \leq C T^{-\frac{\gamma}{2}}$.
Therefore,
\begin{equation}\label{secondstep}
\PP\big[\sup_{t\in[0,T]}\Delta_ {t}\leq 1\big]\leq 
CT^{-\frac{\gamma}{2}}(\|f_1\|_{q_T}+\|f_2\|_{q_T}).
\end{equation}

It remains to prove that $\|f_1\|_{q_T}$ and $\|f_2\|_{q_T}$ are bounded by logarithmic terms. 
$$
\|f_1\|_{q_T}
=a_T^{\alpha_T}  \EE\cro{ V_{a_T}^{q_T/2} \ind_{V_{a_T} \geq 1} + \ind_{V_{a_T} \leq 1}}^{1/q_T}
\leq a_T^{\alpha_T} \pare{\EE\cro{ V_{a_T}^{q_T}}^{1/(2q_T)} + 1}  .$$
By (H2), $V_{a_T}\overset{\mathcal L}{=}a_T^{2-\gamma}V_1$. Therefore,
$$ \EE\cro{ V_{a_T}^{q_T}}^{1/(2q_T)} = a_T^{1-\frac{\gamma}{2}}  \EE\cro{ V_{1}^{q_T}}^{1/(2q_T)} \, .
$$
As in \eqref{momentDelta}, it follows from (H4) that
$$\EE[V_1^{m}] \leq   C c^{-m/\alpha} \frac{m}{\alpha} \Gamma\pare{ \frac{m}{\alpha}} 
\textrm{ for every }m\in\NN,$$
so that using Stirling's formula, it is easy to show that for $T$ large enough
$$\EE[ V_{1}^{q_T}]^{1/2q_T}\leq C(\ln T)^{\frac{1}{2\alpha}}.$$
We conclude that for $T$ large enough
\begin{equation}\label{normef1}
\|f_1\|_{q_T}\leq C(\ln T)^{a(2-\frac{\gamma}{2})+\frac{1}{2\alpha}}.
\end{equation}


Let us now turn our attention to $\|f_2\|_{q_T}$.
$$
\begin{aligned}
\|f_2\|_{q_T} 
& \leq T^{\alpha_T} \EE \cro{e^ {-2q_T\alpha_T V_{a_T}^{1/2}}}^{1/2q_T}
  \EE \cro{\max(V_{a_T},1)^{q_T}}^{1/2q_T}
\\
&   \leq C T   (\ln T)^{a(1-\frac{\gamma}{2})+\frac{1}{2\alpha}} 
 \EE \cro{e^ {-2q_T\alpha_T a_T^{1-\frac \gamma 2} \sqrt{V_1}}}^{1/2q_T}
 \, .
 \end{aligned}
 $$
Let us note $\lambda_T=2q_T\alpha_T a_T^{1-\frac{\gamma}{2}}=2(\ln T)^{1+a(1-\frac{\gamma}{2})}.$
Then using  (H5),
$$\begin{aligned}
\EE \cro{e^ {-\lambda_T V_{1}^{1/2}}} 
& \leq \int_{0}^{+\infty}  \PP \cro{  V_1 \leq \frac{u^2}{\lambda_T^2} }e^{-u} \, du
\\
& \leq C \int_{0}^{+\infty} e^{-c \frac{\lambda_T^{2\beta}}{u^{2\beta}}} e^{-u} \, du \, , 
\end{aligned}$$
for some constants $c$ and $C$. We perform the change of variable $u= \lambda_T^{\frac{2 \beta}{1+2\beta}} v$
in the preceding integral. This yields

$$\begin{aligned}
\EE \cro{e^ {-\lambda_T V_{1}^{1/2}}}
&\leq C \lambda_T^{2\beta/(1+2\beta)} \int_{0}^{+ \infty} e^{-\lambda_T^{\frac{2\beta}{1+2\beta}} 
\pare{v + c v^{-2\beta}}} dv
\\
& \leq C \lambda_T^{2\beta/(1+2\beta)}\int_{1}^{+ \infty}  e^{-\lambda_T^{\frac{2\beta}{1+2\beta}} v} \, dv 
+ C \lambda_T^{2\beta/(1+2\beta)} \int_{0}^{1} e^{-\lambda_T^{\frac{2\beta}{1+2\beta}} c v^{-2\beta}}
\, dv 
\\
& \leq C  e^{-\lambda_T^{\frac{2\beta}{1+2\beta}}} +  C  e^{-c \lambda_T^{\frac{2\beta}{1+2\beta}}}  
\\
& \leq C e^{-c \lambda_T^{\frac{2\beta}{1+2\beta}}}  \, , 
 \end{aligned}$$
for some other constants $c$ and $C$. This leads to 
$$\|f_2\|_{q_T} \leq CT(\ln T)^{a(1-\frac{\gamma}{2})+\frac{1}{2\alpha}} 
e^{-c (\ln T)^{\frac{2\beta}{1+2\beta} \pare{1+a (1-\frac{\gamma}{2})}-1}} \, .$$
We choose $a$ such that $a(1-\frac{\gamma}{2}) > 1 + \frac{1}{\beta}$.  Then
 $T e^{-c (\ln T)^{\frac{2\beta}{1+2\beta} \pare{1+a (1-\frac{\gamma}{2})}-1}}$ remains bounded, and 
when $T$ goes to infinity, we get
\begin{equation}\label{normef2}
\|f_2\|_{q_T}\leq C(\ln T)^c \, ,
\end{equation}
for some constant $c >1 + \frac{1}{\beta}+\frac{1}{2\alpha}$. 
From \eqref{secondstep}, \eqref{normef1} and \eqref{normef2} we deduce that for $T$ large enough,
$$\PP \big[\sup_{t\in[0,T]}\Delta_ {t}\leq 1\big]
\leq CT^{-\frac{\gamma}{2}}(\ln T)^{c},$$
with $c> \frac{\beta+1}{\beta} \frac{4-\gamma}{2-\gamma} +\frac{1}{2\alpha}$.

     \medskip
    
{\noindent\bf Acknowledgments.} We are grateful to Yueyun Hu for stimulating discussions.


%

\end{document}